%
\documentstyle{amsppt}
\magnification=1200

\baselineskip=16.5pt
\tolerance=5000
\null
\vskip2truecm
\def\irtp{\Cal P}
\def\irta{\Cal A}

\topmatter
\title On the cardinality and weight spectra of compact spaces, II\endtitle
\author I. Juh\'asz$^*$ and S. Shelah$^{**}$\endauthor
\thanks $^*$Research supported by the Hungarian National Foundation
for Scientific Research grant no. 16391.\newline\indent
$^{**}$Research supported by ``The Israel Science Foundation'',
administered by The Israel Academy of Sciences and Humanities, Publication
no\. 612\endthanks
\abstract Let $B(\kappa,\lambda)$ be the subalgebra of $\irtp(\kappa)$
generated by $[\kappa]^{\le\lambda}$. It is shown that if $B$ is any
homomorphic image of $B(\kappa,\lambda)$ then either $|B|<2^\lambda$ or
$|B|=|B|^\lambda$, moreover if $X$ is the Stone space of $B$ then either
$|X|\le 2^{2^\lambda}$ or $|X|=|B|=|B|^\lambda$.\newline
\indent This implies the existence of 0-dimensional compact $T_2$ spaces whose
cardinality and weight spectra omit lots of singular cardinals of ``small''
cofinality.\endabstract
\endtopmatter

\document
\baselineskip=16.5pt

\noindent{\bf 1. Introduction}

\smallskip

It was shown in [J] that for every uncountable regular cardinal $\kappa$, if
$X$ is any compact $T_2$ space with $w(X)>\kappa$ ($|X|>\kappa$) then $X$ has a
closed subspace $F$ such that $\kappa\le w(F)\le 2^{<\kappa}$
$\left(\text{resp. }
\kappa\le|F|\le\sum\{2^{2^\lambda}\colon\lambda<\kappa\}\right)$.
In particular, the weight or cardinality spectrum of a compact space may never
omit an inaccessile cardinal, moreover under GCH the weight spectrum cannot
omit any uncountable regular cardinal at all.

In the present note we prove a theorem which implies that for singular $\kappa$
on the other hand there is always a 0-dimensional compact $T_2$ space whose
cardinality and weight spectra both omit $\kappa$.

We formulate our main result in a boolean algebraic framework. The topological
consequences easily follow by passing to the Stone spaces of the boolean
algebras that we construct.

\bigskip

\noindent{\bf 2. The Main Result}

We start with a general combinatorial lemma on binary relations. In order to
formulate it, however, we need the following definitions.

\definition{Definition 1} Let $\prec$ be an arbitrary binary relation on a set
$X$ and $\tau,\mu$ be cardinal numbers. We say that $\prec$ is $\tau$-full if
for every subset $a\subset X$ with $|a|=\tau$ there is some $x\in X$ such that
$|\{y\in a\colon y\prec x\}|=\tau$. Moreover, $\prec$ is said to be $\mu$-local
if for every $x\in X$ we have $|\text{pred}\,(x,\prec)|\leqq\mu$, where
$\text{pred}(x,\prec)=\{y\in X\colon y\prec x\}$.
\enddefinition

Now, our lemma is as follows.

\proclaim{Lemma 2} Let $\prec$ be a binary relation on the cardinal $\varrho$
that is both $\tau$-full and $\mu$-local. Then for every almost disjoint family
$\irta\subset[\varrho]^\tau$ we have
$$
\vert\irta\vert\le\varrho\cdot\mu^\tau.
$$
\endproclaim

\demo{Proof} For every set $a\in\irta$ there is a $\xi_a\in\varrho$ such that
$g(a)=a\cap\text{pred}(\xi_a,\prec )$ has cardinality $\tau$ because $\prec$ is
$\tau$-full. This map $g$ is clearly one-to-one for $\irta$ is almost disjoint.
But the range of $g$ is a subset of
$\cup\{[\text{pred}(\xi,\prec)]^\tau\colon\xi\in\varrho\}$ whose cardinality
does not exceed $\varrho\cdot[\mu]^\tau$, and this completes the proof.

Before we formulate our main result we need some notation. Given the cardinals
$\kappa$ and $\lambda$ (we may assume $\lambda\le\kappa$) we denote by
$B(\kappa,\lambda)$ the boolean subalgebra of the power set algebra
$\irtp(\kappa)$ generated by all subsets of $\kappa$ os size $\le\lambda$. In
other words
$$
B(\kappa,\lambda)=[\kappa]^{\le\lambda}\cup
\{x\subset\kappa\colon\kappa\setminus x\in[\kappa]^{\le\lambda}\}.
$$
What we can show is that the size of a homomorphic image of $B(\kappa,\lambda)$
(as well as the size of its Stone space) has to satisfy certain restrictions,
namely it is either ``small'' or cannot have ``very small'' cofinality.
\enddemo

\proclaim{Theorem 3} Let $h\colon B(\kappa,\lambda)\to B$ be a homomorphism of
$B(\kappa,\lambda)$ onto the boolean algebra $B$. Then (i) either
$|B|<2^\lambda$ or $|B|^\lambda=|B|$; (ii) if $X=St(B)$ is the Stone space of
$B$ then either $|X|\le 2^{2^\lambda}$ or $|X|=|B|=|B|^\lambda$.
\endproclaim

\demo{Proof} Set $|B|=\varrho$ and assume that $\varrho\ge 2^\lambda$. Since
$[\kappa]^{\le\lambda}$ generates $B(\kappa,\lambda)$ therefore
$A=h''[\kappa]^{\le\lambda}$ generates $B$ and thus we have $|A|=\varrho$ as
well. We claim that the relation $\le_B$ is
\roster
\item"{(a)}" $\tau$-full on $A$ for each $\tau\le\lambda$;
\item"{(b)}" $2^\lambda$-local on $A$.
\endroster

Indeed, if $a\in[A]^\tau$ where $\tau\le\lambda$ then there is a set
$x\in[[\kappa]^{\le\lambda}]^\tau$ such that $a=h''x$. But then $b=\cup
x\in[\kappa]^{\le\lambda}$ as well, hence $h(b)\in A$ and clearly $a\subset
\text{pred}(h(b),\le_B)$ because $h$ is a homomorphism. This, of course, is
much more than what we need for (a).

To see (b), let us first note that if $b,c\in[\kappa]^{\le\lambda}$ and
$h(b)\le h(c)$ then $b\cap c\in[\kappa]^{\le\lambda}$ as well and $h(b\cap
c)=h(b)\land h(c)=h(b)$ using that $h$ is a homomorphism again. But this
implies $\text{pred} (h(c),\le_B)=h''\irtp(c)$ for any 
$c\in[\kappa]^{\le\lambda}$,
consequently $|\text{pred}(h(c),\le_B)|\le|\irtp(c)|\le 2^\lambda$
and this completes the proof of (b).

Applying Lemma 2 we may now conclude that for every cardinal $\tau\le\lambda$
and for every almost disjoint family $\irta\subset[\varrho]^\tau$ we have
$$
|\irta|\le\varrho\cdot (2^\lambda)^\tau=\varrho.
$$
This, in turn, implies $\varrho^\lambda=\varrho$. Indeed, assume that
$\varrho^\lambda>\varrho$ and $\tau $ be the smallest cardinal with
$\varrho^\tau>\varrho$. Then $\tau\le\lambda$ and $\varrho^{<\tau}=\varrho$,
and as is well-known, there is an almost disjoint family
$\irta\subset[{}^{<\tau}\varrho]^\tau$ of size $\varrho^\tau>\varrho$, namely
$\irta=\{A_f\colon f\in {}^\tau\varrho\}$ where $A_f=
\{f\restriction\xi\colon\xi<\tau\}$ for any $f\in {}^\tau\varrho$.

Now, to prove (ii) first note that if $|B|\le 2^\lambda$ then trivially $|X|\le
2^{2^\lambda}$. So assume $|B|>2^\lambda$ and in this case we prove that
actually
$$
|X|=2^{2^\lambda}\cdot|B|.
$$
We first show that $|X|\geqq 2^{2^\lambda}\cdot |B|$, which, as $|X|\ge|B|$ is
always valid, boils down to showing that $|X|\ge 2^{2^\lambda}$.

Using that $|B|=|h''[\kappa]^{\le\lambda}|=\varrho>2^\lambda$
we may select a collection $\{a_\alpha\colon\alpha\in(2^\lambda)^+\}
\subset[\kappa]^{\le\lambda}$ such that $\alpha\ne \beta$ implies
$h(a_\alpha)\ne h(a_\beta)$ and by a straight forward $\Delta$-system argument
we may also assume that $\{a_\alpha\colon\alpha\in(2^\lambda)^+\}$
is a $\Delta$-system with root $a$. Then, as $h$ is a homomorphism, we also
have $h(a_\alpha)\land h(a_\beta)=h(a)$ for distinct $\alpha$ and $\beta$ and
so $\{h(a_\alpha)-h(a)\colon\alpha\in(2^\lambda)^+\}$ are pairwise disjoint and
disdinct elements $B$, all but at most one of which is non-zero. However the
existence of $2^\lambda$ many pairwise disjoint non-zero elements in a boolean
algebra clearly implies the existence of $2^{2^\lambda}$ ultrafilters in it,
hence we are done with showing $|X|\ge 2^{2^\lambda}$.

Next, to see $|X|\le 2^{2^\lambda}\cdot |B|$ note that, again as $h$ is a
homomorphism, $h''[\kappa]^{\le\lambda}$ is a (not necessarily proper) ideal
in $B$, hence there is no more than one ultrafilter $u$ on $B$ such that $u\cap
h''[\kappa]^{\le\lambda}=\emptyset$. If, on the other hand, $u\in X$ is such
that $b\in u\cap h^{cc[\kappa]^{\le\lambda}}$ then $u$ is generated by its
subset $u\cap\text{pred}(b,\le_B)$. However $\le_B$ is clearly $2^\lambda$-
local on $h''[\kappa]^{\le\lambda}$, and so we conclude that
$$\aligned
|X|\le & 1+|\cup
\{\irtp(\text{pred}(b,\le_B))\colon b\in h''[\kappa]^{\le\lambda}\}|\le\\
&\le 1+2^{2^\lambda}\cdot|B|=2^{2^\lambda}\cdot|B|.\endaligned
$$
This completes the proof of our theorem.

Now let $X(\kappa,\lambda)$ be the Stone space of the boolean algebra
$B(\kappa,\lambda)$. Using Stone duality and the notation of [J] the above
result has the following reformulation about the weight and cardinality spectra
of the 0-dimensional compact $T_2$ space $X(\kappa,\lambda)$.
\enddemo

\proclaim{Corollary 4}
\roster
\item"{(i)}" For every $\mu\in Sp(w,X(\kappa,\lambda))$ we have either
$\mu<2^\lambda$ or $\mu^\lambda=\mu$, hence $cf(\mu)>\lambda$;
\item"{(ii)}" if $\mu\in Sp(|\,\,|,X(\kappa,\lambda))$ then either
$\mu<2^{2^\lambda}$ or $\mu^\lambda=\mu$.
\endroster
In fact, for every closed subspace $Y$ of $X(\kappa,\lambda)$ we have either
$w(Y)\le 2^\lambda$ or $w(Y)^\lambda=w(Y)$ and $|Y|=2^{2^\lambda}\cdot w(Y)$.
\endproclaim

It follows from this immediately that if $2^{2^\lambda}<\kappa$ then the
cardinality and weight spectra of the space $X(\kappa,\lambda)$ omit every
cardinal $\mu\in(2^{2^\lambda},\kappa]$ with $cf(\mu)\le\lambda$. In
particular, if GCH holds then $\lambda<\kappa$ implies that both $Sp(|\,\,|,
X(\kappa,\lambda))$ and $Sp(w,X(\kappa,\lambda))$ omit all cardinals
$\mu\in(\lambda,\kappa]$ with $cf(\mu)\le\lambda$.

Note that similar omission results were obtained by van Douwen in [vD] for the
case $\lambda=\omega$ and ${\kappa}$ strong limit.

An interesting open problem arises here that we could not settle: Can one find
for every cardinal $\kappa$ a compact $T_2$ space $X$ such that the cardinality
and/or weight spectra of $X$ omit every singular cardinal below $\kappa$?

\vskip2truecm

\centerline{{\smc References}}

\smallskip

\roster
\item"{[vD]}" E\. van Douwen, {\it Cardinal functions
on compact F-spaces and on weakly countably compact boolean algebras,}
Fund\. Math. {\bf 114} (1981), 236-256.
\item"{[J]}" I\. Juh\'asz, {\it On the weight spectrum of a compact
spaces,} 
Israel J\.
Math. {\bf 81} (1993), 369-379.
\endroster
\enddocument